\documentclass[11pt]{article}

\usepackage{amsmath, amsthm, amssymb}
\usepackage{verbatim}

\date{}
\title{\vspace{-0.8cm}Resilient pancyclicity of random and pseudo-random graphs}
\author{
Michael Krivelevich \thanks{School of Mathematical Sciences, Raymond
and Beverly Sackler Faculty of Exact Sciences, Tel Aviv University,
Tel Aviv, 69978, Israel. Email: krivelev@post.tau.ac.il. Research
supported in part by USA-Israel BSF grant 2006322, by grant 1063/08
from the Israel Science Foundation, and by a Pazy memorial award.}
\and Choongbum Lee \thanks{Department of Mathematics, UCLA, Los
Angeles, CA, 90095. Email: choongbum.lee@gmail.com. Research
supported in part by Samsung Scholarship.} \and Benny Sudakov
\thanks{Department of Mathematics, UCLA, Los Angeles, CA 90095.
Email: bsudakov@math.ucla.edu. Research supported in part by NSF
CAREER award DMS-0812005 and by a USA-Israel BSF grant. } }

\theoremstyle{plain}
\newtheorem{thm}{Theorem}[section]
\newtheorem{prop}[thm]{Proposition}
\newtheorem{lemma}[thm]{Lemma}

\newtheorem{conj}[thm]{Conjecture}

\theoremstyle{definition}
\newtheorem{dfn}{Definition}

\newenvironment{pf}{\noindent\textbf{Proof.}}{\qed\medskip}
\newenvironment{rem}{\noindent\textbf{Remark.}}{\medskip}

\begin{document}
\maketitle

\begin{abstract}
A graph $G$ on $n$ vertices is \textit{pancyclic} if it contains
cycles of length $t$ for all $3 \leq t \leq n$. In this paper we
prove that for any fixed $\epsilon>0$, the random graph $G(n,p)$
with $p(n)\gg n^{-1/2}$ asymptotically almost surely has the
following resilience property. If $H$ is a subgraph of $G$ with
maximum degree at most $(1/2 - \epsilon)np$ then $G-H$ is pancyclic.
In fact, we prove a more general result which says that if $p \gg
n^{-1+1/(l-1)}$ for some integer $l \geq 3$ then for any
$\epsilon>0$, asymptotically almost surely every subgraph of
$G(n,p)$ with minimum degree greater than $(1/2+\epsilon)np$
contains cycles of length $t$ for all $l \leq t \leq n$. These
results are tight in two ways. First, the condition on $p$
essentially cannot be relaxed. Second, it is impossible to improve
the constant $1/2$ in the assumption for the minimum degree. We also
prove corresponding results for pseudo-random graphs.

\end{abstract}

\section{Introduction}
\label{introduction_section}

A typical result in graph theory can be stated as ``Under certain
conditions, a graph $G$ possesses a property $\mathcal{P}$''. Once
this type of a result is established, it is natural to ask ``How
strongly does $G$ possess $\mathcal{P}$?''. In fact, several
important results in extremal graph theory can be viewed as an
answer to this question for various graph properties (we will
provide some concrete examples after introducing necessary
definitions). In this paper we will study this question in the
context of random and pseudo-random graphs. The random graph model
we consider is the binomial random graph $G(n,p)$. The random graph
$G(n,p)$ denotes the probability space whose points are graphs with
vertex set $[n] = \{1,\ldots,n\}$ where each pair of vertices forms
an edge randomly and independently with probability $p$. We say that
$G(n,p)$ possesses a graph property $\mathcal{P}$
\textit{asymptotically almost surely}, or a.a.s. for brevity, if the
probability that $G(n,p)$ possesses $\mathcal{P}$ tends to 1 as $n$
goes to infinity. The pseudo-random graphs we will study are
$(n,d,\lambda)$-graphs with $\lambda = o(d)$, where an
$(n,d,\lambda)$-graph is a $d$-regular graph on $n$ vertices whose
second largest (in absolute value) eigenvalue of the adjacency
matrix is bounded by $\lambda$. The abundance of structure and
results arising from this simple looking definition is quite
surprising (see, e.g., \cite{MR2223394} for more details). A graph
property is called \textit{monotone increasing (decreasing)} if it
is preserved under edge addition (deletion).

The main concept studied in this paper and briefly outlined above is
that of \textit{resilience}. Formally, following \cite{MR2462249},
we define:

\begin{dfn} Let $\mathcal{P}$ be a monotone increasing (decreasing) graph property.
\begin{itemize}
\item[(i)] (Global resilience) The global resilience of $G$ with respect to $\mathcal{P}$ is the minimum number $r$ such that by deleting (adding) $r$ edges from $G$ one can obtain a graph not having $\mathcal{P}$.
\item[(ii)] (Local resilience) The local resilience of a graph $G$ with respect to $\mathcal{P}$ is the minimum number $r$ such that by deleting (adding) at most $r$ edges at each vertex of $G$ one can obtain a graph not having $\mathcal{P}$.
\end{itemize}
\end{dfn}

Using this terminology, one can state the celebrated theorem of
Tur\'{a}n \cite{MR0018405} as ``The complete graph on $n$ vertices
$K_n$ has global resilience $\frac{n^2}{2r} - \frac{n}{2}$ with
respect to being $K_{r+1}$-free''. Another classical theorem, that
of Dirac (see, e.g., \cite{MR2159259}) can be rephrased as ``$K_n$
has local resilience $\left\lfloor n/2 \right\rfloor$ with respect
to Hamiltonicity''. As these examples suggest, the notion of
resilience lies in the center of extremal graph theory. In
\cite{MR2462249}, Sudakov and Vu have initiated the systematic study
of global and local resilience of random and pseudo-random graphs.
They obtained resilience results with respect to various properties
such as perfect matching, hamiltonicity, chromatic number and having
a nontrivial automorphism (this result appeared in their earlier
paper with Kim \cite{MR1945368}). For example, they showed that if
$p > \log^4 n / n$ then a.a.s. any subgraph of $G(n,p)$ with minimum
degree $(1/2+o(1))np$ is hamiltonian. An interesting thing to notice
is that this result can be viewed as a generalization of Dirac's
Theorem mentioned above as a complete graph is also a random graph
$G(n,p)$ with $p=1$. As we will see, this connection is very natural
and most of the resilience results can be viewed as a generalization
of classic graph theory results to random and pseudo-random graphs.

There are several other papers that obtained resilience type
results. Krivelevich and Frieze \cite{MR2430433} gave a lower bound
(not tight) on resilience of $G(n,p)$ with respect to  being
Hamiltonian in the range of $p$ not covered by the above mentioned
result of Sudakov and Vu. Dellamonica, Kohayakawa, Marciniszyn and
Steger \cite{MR2383452} studied the global resilience of random
graphs with respect to containing a cycle of length at least
$(1-\alpha)n$ for a fixed $\alpha$ as a generalization of a theorem
of Woodall \cite{MR0318000}.  Recently, Ben-Shimon, Krivelevich and
Sudakov \cite{MR000000} investigated the resilience of random
regular graphs with respect to being hamiltonian.

A graph on $n$ vertices is called \textit{pancyclic} if it contains
cycles of length $t$ for all $3 \leq t \leq n$. In this paper, we
study the resilience of random and pseudo-random graphs with respect
to this property. Similarly to the above mentioned results, our
result can also be viewed as a generalization of a classical result
in graph theory -- that by Bondy \cite{MR0285424}. It says that if
$G$ is a graph on $n$ vertices with minimum degree greater than
$n/2$, then $G$ is pancyclic. The corresponding theorems we prove
are,

\begin{thm} \label{pancyclic_thm3}
If $p \gg n^{-1/2}$ then $G(n,p)$ asymptotically almost surely has
local resilience $(1/2 + o(1))np$ with respect to being pancyclic.
\end{thm}
\begin{thm} \label{pancyclic_pseudo_thm00}
Let $G=(V,E)$ be a $(n,d,\lambda)$-graph satisfying $d^2 /n \gg \lambda$. Then $G$ has local resilience $(1/2 + o(1))d$ with respect to being pancyclic.
\end{thm}

Our results are asymptotically tight in two ways. First, one cannot
improve the constant $1/2$ since both random and pseudo-random
graphs can be made bipartite by randomly partitioning the graph into
two equal size parts. In this way we typically have a subgraph with
minimum degree about one half of the original degree which does not
contain any odd cycles. Second, the restrictions on the parameters
are also essentially tight. To see this for random graphs, note that
if $p \ll n^{-1/2}$ then typically each vertex has degree
$(1+o(1))np$ and the number of triangles containing each vertex is
at most ${O}(n^2p^3) \ll np$. Therefore deleting edges of all
triangles leaves all degrees essentially unchanged. For
pseudo-random graphs this can be derived from a variant of the
construction of Alon \cite{MR1302331} (see, e.g., \cite{MR2223394})
which gives a triangle-free $(n,d,\lambda)$-graph with $d^3 /
\lambda^2 = \Theta (n)$.

We can also prove more general results for sparser graphs. Let the
\textit{girth} of a graph be the length of its shortest cycle and
the \textit{circumference} be the length of its longest cycle.
Brandt, Faudree and Goddard \cite{MR1611825} called a graph
\textit{weakly pancyclic} if it contains cycles of length $t$ where
$t$ ranges from its girth up to its circumference. The following
theorems, generalizing Theorems \ref{pancyclic_thm3} and
\ref{pancyclic_pseudo_thm00}, are motivated by this concept of weak
pancyclicity.

\begin{thm} \label{pancyclic_thm4}
 For any fixed integer $l \geq 3$, if $p \gg n^{-1+1/(l-1)}$ then $G(n,p)$
 asymptotically almost surely has local resilience $(1/2 + o(1))np$ with respect to containing cycles of length $t$ for all $l \leq t \leq n$.
\end{thm}
\begin{thm} \label{pancyclic_pseudo_thm01}
Let $k$ be either 3 or an even integer satisfying $k \geq 4$ and let $G=(V,E)$ be a $(n,d,\lambda)$-graph satisfying $d^{k-1} /n \gg \lambda^{k-2}$. Then  $G$ has local resilience $(1/2 + o(1))d$ with respect to containing cycles of length $t$ for all $k \leq t \leq n$.
\end{thm}

The above results are not exactly weak pancyclicity results since if
we allow the adversary to delete half of the edges at each vertex he
might decide not to remove a 3-cycle and then remove every other
cycle of length $4$ up to $l-1$. But still it is best to view these
results in the context of weak pancyclicity. Similarly as before,
the result for random graphs is asymptotically tight. Indeed, note
that if $p \ll n^{-1+1/(l-1)}$ then typically each vertex of the
random graph has degree $(1+o(1))np$ and the number of cycles of
length $l$ containing each vertex is at most $O(n^{l-1}p^l) \ll np$.
Therefore we can delete few edges from each vertex to remove every
$l$-cycle. We suspect that our result for pseudo-random graphs is
asymptotically tight as well. Note that the assumption $d^{k-1} /n
\gg \lambda^{k-2}$ in particular implies $\lambda = o(d)$ since
$d^{k-2} \geq d^{k-1} /n \gg \lambda^{k-2}$, so even when we do not
explicitly mention $\lambda = o(d)$, we are always in this
situation. Although odd integers $k > 3$ are omitted from the result
of pseudo-random graphs, nevertheless in this case $d^{k-1} /n \gg
\lambda^{k-2}$ implies $d^{k} /n \gg \lambda^{k-1}$ and so by using
the result for $k+1$ (which is now even) we can find cycles of
length $t$ for all $k+1 \leq t \leq n$. We believe that the result
of Theorem \ref{pancyclic_pseudo_thm01} is valid also for odd $k\ge
5$, but at present we do not have enough tools to verify it. We will
address this point in more details in concluding remarks.

The rest of this paper is organized as follows. In Section
\ref{preliminaries_section} we collect some known results which we
need later to prove our main theorems. In Section
\ref{randomgraphproperty_section} we establish properties of random
graphs and use them in Section \ref{proofofrandomgraphthm_section}
to prove Theorem \ref{pancyclic_thm4}. In Sections
\ref{pseudorandomgraphproperty_section},
\ref{proofofpseudorandomgraphthm_section} we follow the same pattern
to prove the pseudo-random graph analog, Theorem
\ref{pancyclic_pseudo_thm01}. The last section contains some
concluding remarks and open problems.
\\

\noindent\textbf{Notation.} $G=(V,E)$ denotes a graph with vertex
set $V$ and edge set $E$. We use $v \sim w$ to indicate that $v,w$
are adjacent. $\Delta(G), \delta(G)$ denote the maximum degree and
the minimum degree of $G$, respectively. For a set $X \subset V$,
let $N(X)$ be the collection of all vertices $v$ which are adjacent
to at least one vertex in $X$. If $X=\{u\}$ is a singleton set we
denote its neighborhood by $N(u)$. Let $N^{(0)}(v) := \{ v \}$ and
$N^{(k)}(v)$ be the vertices at distance exactly $k$ from $v$. This
can also be recursively defined as $N^{(k)}(v) = N(N^{(k-1)}(v))
\backslash (N^{(k-1)}(v) \cup N^{(k-2)}(v))$. Note that $N^{(1)}(v)
= N(v)$. For a set $X$, we denote by $E(X)$ the set of edges in the
induced subgraph $G[X]$ and by $e(X)= |E(X)|$ its size. Similarly,
for two sets $X$ and $Y$, we denote by $E(X,Y)$ the set of ordered
pairs $(x,y) \in E$ such that $x \in X$ and $y \in Y$, also $e(X,Y)
= |E(X,Y)|$. Note that $e(X,X) = 2e(X)$. If we have several graphs,
then the graph we are currently working with will be stated as a
subscript. For example $N^{(k)}_G(v)$ is the $k$-th neighborhood of
$v$ in graph $G$. A cycle of length $l$ is denoted by $C_l$.

We also utilize the following standard asymptotic notation. For two
functions $f(n)$ and $g(n)$, write $f(n)=\Omega(g(n))$ if there
exists a constant $C$ such that $\liminf_{n \rightarrow \infty}
f(n)/g(n) \geq C$. If there is a subscript such as in
$\Omega_\epsilon$ this means that the constant $C$ may depend on
$\epsilon$. We write $f(n)=o(g(n))$ or $f(n) \ll g(n)$  if
$\limsup_{n\rightarrow \infty}f(n)/g(n) = 0$. Also, $f(n)=O(g(n))$
if there exists a positive constant $C>0$ such that
$\limsup_{n\rightarrow \infty} f(n)/g(n) \leq C$. Throughout the
paper log denotes the natural logarithm. To simplify the
presentation, we often omit floor and ceiling signs whenever these
are not crucial and make no attempts to optimize absolute constants
involved. We also assume that the order $n$ of all graphs tends to
infinity and therefore is sufficiently large whenever necessary.

\section{Preliminaries}
\label{preliminaries_section}

In this section we collect various results to be used later in the
proofs of the theorems.

\subsection{Resilience}

The local resilience of random graphs with respect to being
hamiltonian \cite{MR2462249} and containing fixed cycles
(\cite{MR1339852}, \cite{MR1394514}, \cite{MR2145507}) have been
studied before and our arguments for the proof of the main theorems
will use these results. The following results about the local
resilience of random and pseudo-random graphs with respect to
hamiltonicity were proved in \cite{MR2462249}.

\begin{thm} \label{pancyclic_thm1}
 For every fixed $\epsilon > 0$, if $p \geq \log^4n /n$ then the random graph $G(n,p)$ with probability $1 - o(n^{-1})$ has local resilience at least $(1/2 - \epsilon)np$ with respect to being hamiltonian.
\end{thm}
\begin{rem} The above formulation is stronger than the original statement since it explicitly
states the success probability to be $1 - o(n^{-1})$. But this
conclusion follows from the original argument if one carefully
performs the error probability calculations. We will need this
stronger estimate on success probability for our application.
\end{rem}

During the proof we will work with graphs that are similar to
$(n,d,\lambda)$-graphs but are not necessarily regular. The
particular graphs we will encounter are graphs $G=(V,E)$ on $n$
vertices that have minimum degree at least $(1 - \epsilon)d$ and
satisfy the constraint
$$\big|e(X,Y) - \frac{d}{n}|X||Y|\big| \leq
\lambda \sqrt{|X||Y|} \mbox{ for all } X,Y \subset V
$$
on the number of edges between sets. We will call such graphs
\textit{$(n,\epsilon, d,\lambda)$-graphs}. Observe that
$(n,0,d,\lambda)$-graphs are a more general/flexible concept than
that of $(n,d,\lambda)$-graphs, as it does not put specific
assumptions on graph eigenvalues.

\begin{thm} \label{pancyclic_pseudo_thm2} Fix $\epsilon, \epsilon'$ such that $0 \leq 5\epsilon' < \epsilon, k \geq 3$ and let $G$ be an $(n,\epsilon', d,\lambda)$-graph satisfying $d^{k-1}/n = \omega(n) \lambda^{k-2}$ for an arbitrary function $\omega(n)$ increasing to infinity. Then for large enough $n$, $G$ has local resilience at least $(1/2 - \epsilon)d$ with respect to being hamiltonian.
\end{thm}
\begin{rem}
The theorem above does not appear in the original paper \cite{MR2462249} and unfortunately we cannot directly apply the result from Sudakov and Vu. But in fact, they proved a general theorem which can be modified to work under the assumption above. The necessarily modification will be given in the Appendix.
\end{rem}

Next we state the results of Haxell, Kohayakawa and {\L}uczak (\cite{MR1339852}, \cite{MR1394514})
 about the local resilience of random graphs with respect to containing a fixed cycle $C_l$, and of Sudakov, Szab\`{o} and Vu \cite{MR2145507} about the local resilience of pseudo-random graphs with respect to containing a triangle.

\begin{thm} \label{pancyclic_thm2}
For any fixed integer $l \geq 3$ and $\epsilon > 0$, there exists a constant $C=C(l, \epsilon)$ such that, if $p \geq C n^{-1+1/(l-1)}$ then $G(n,p)$ a.a.s. has local resilience at least $(1/2 - \epsilon)np$ with respect to containing $C_l$.
\end{thm}

\begin{thm} \label{pancyclic_pseudo_thm1}
Let $G$ be a $(n,d,\lambda)$-graph satisfying $d^2 / n \geq
\omega(n) \lambda$ for an arbitrary function $\omega(n)$ tending to
infinity. Then $G$ has local resilience $(1/2 + o(1))d$ with respect
to containing a triangle.
\end{thm}
\begin{rem}
Both theorems are originally stated in a global resilience form but
for convenience we stated it as above in a slightly weaker local
resilience form. Also the conclusion of Theorem \ref{pancyclic_thm2}
(as stated) for even cycles is weaker than in the original paper.
\end{rem}

\subsection{Extremal Graph Theory}

The following simple but useful lemma allows one to find a large
minimum degree subgraph in a graph with large average degree. (See,
e.g., \cite{MR2159259}, Proposition 1.2.2)

\begin{lemma} \label{lemma_findingmindegree}
Let $G=(V,E)$ be a graph on $n$ vertices with at least $dn/2$ edges. Then $G$ contains a subgraph $G' \subset G$ with minimum degree at least $d/2$.
\end{lemma}

Next theorem is a classical result by Bondy and Simonovits \cite{MR0340095} about even cycles in graphs.
\begin{thm} \label{bondysimonovits1_thm}
Let $k$ be a positive integer and $G=(V,E)$ be a graph on $n$ vertices satisfying $|E| > 90kn^{1 + 1/k}$. Then $G$ contains a cycle of length $2k$.
\end{thm}

We will also need the celebrated P\`{o}sa rotation-extension lemma
(see \cite{MR2321240}, Ch. 10, Problem 20). This lemma will help us
in finding long paths in a graph with  expansion properties.
\begin{lemma} \label{lemma_pancyclic3}
Let $G=(V,E)$ be a graph
such that $|N(X)\setminus X| \geq 2|X| -1$ for all $X \subset V$
with $|X| \leq t$. Then for any vertex $v \in V$ there exists a path
of length $3t-2$ in $G$ that has $v$ as an end point.
\end{lemma}

\subsection{Concentration}

The following two well-known concentration results (see, for example
\cite{MR1782847}, Theorems 2.3 and 2.10) will be used several times
during the proof.
We denote by $Bi(n,p)$ a binomial random variable with parameters
$n$ and $p$.
\begin{thm}(Chernoff inequality) If $X \sim Bi(n,p)$ and $\epsilon>0$, then
\[ P\big(|X - \mathbb{E}[X]| \geq \epsilon \mathbb{E}[X]\big) \leq e^{-\Omega_\epsilon(\mathbb{E}[X])}. \]
\end{thm}

Let $m,n$ and $N$ be positive integers with $m,n < N$, let $X = [N],
X' = [n]$, and let $A$ be a $m$-element subset of $X$ chosen
uniformly at random. Then the distribution of the random variable
$|A \cap X'|$ is called the \textit{hypergeometric distribution}
with parameters $N, n$ and $m$.
\begin{thm} \label{lemma_hypergeometric}
Let $X$ have the hypergeometric distribution with parameters $N, n$ and $m$. Then,
\[ P\big(|X - \mathbb{E}[X]| \geq \epsilon \mathbb{E}[X]\big) \leq e^{-\Omega_\epsilon(\mathbb{E}[X])}. \]
\end{thm}


\section{Properties of Random Graphs}
\label{randomgraphproperty_section}

In this section we establish properties of random graphs to be used
later to prove Theorem \ref{pancyclic_thm4}.

First we show formally a rather expected monotonicity property --
(relative) local resilience with respect to cycles can only grow
with the edge probability $p(n)$.

\begin{prop} \label{probrestriction_prop}
Let $l$ be fixed and let $p'=p'(n)$ satisfy: $0 < p' \leq p \leq 1$
and $np' \gg \log n$. If $G(n, p')$ a.a.s. has local resilience at
least $(1/2 - \epsilon/2)np'$ with respect to containing cycles of
length $t$ for all $l \leq t \leq n$ then $G(n,p)$ a.a.s. has local
resilience at least $(1/2 - \epsilon)np$ with respect to the same
property.
\end{prop}
\begin{pf}
Let $\mathcal{P}$ be the property of having local resilience at
least $(1/2 - \epsilon/2)np'$ with respect to containing every
cycles of length $t$ for all $l \leq t \leq n$. Define $q = p' / p$
and consider the following two round process of exposing the edges
of $G(n,p')$. In the first round, every edge appears with
probability $p$ (call this graph $G_1$). Then at the second round,
every edge that appeared in the first round will remain with
probability $q$ and will be deleted with probability $1-q$ (call
this graph $G_2$). Then $G_1$ has the same distribution as $G(n,p)$
and $G_2$ has the same distribution as $G(n,p')$. By our assumption
we know that $G_2$ a.a.s. has property $\mathcal{P}$. Now define $X$
to be the event that $G_1$ satisfies: $P(G_2 \notin \mathcal{P} |
G_1) \geq 1/2$. Then $(1/2)P(X) \leq P(G_2 \notin \mathcal{P})=o(1)$
and therefore $P(X) = o(1)$. Thus a.a.s. in $G(n,p)$, $P(G_2 \notin
\mathcal{P} | G_1) < 1/2$ or in other words $P(G_2 \in \mathcal{P} |
G_1) \geq 1/2$. Let $\mathcal{A}$ be the collection of graphs
$G_1\in G(n,p)$ having this property.

Now given any subgraph $H$ of $G_1$ with maximum degree at most
$(1/2 - \epsilon)np$, select every edge with probability $q$ to get
a graph $H'$. Then by Chernoff inequality, each vertex of $H'$ has
maximum degree at most $(1/2 - \epsilon/2)np'$ with probability at
least $1 - e^{-\Omega_\epsilon(np')} = 1 - o(n^{-1})$. Therefore
$H'$ has maximum degree at most $(1/2 - \epsilon/2)np'$ with
probability at least $1-o(1)$.

 Finally to put things together, condition on the event that $G_1=G(n,p) \in \mathcal{A}$. By the first part of the proof a.a.s. $G \in \mathcal{A}$ so if we can prove the claim under this assumption then we are done. Given a subgraph $H \subset G$ with maximum degree at most $(1/2 - \epsilon)np$, sample every edge of $G$ with probability $q$ to obtain subgraphs $H' \subset G' \subset G$. Since $G \in \mathcal{A}$, we know that $P(G' \in \mathcal{P} | G) \geq 1/2$ and by the second part of the proof we know that $P\big( \Delta(H') \leq (1/2 - \epsilon/2)np'\big) \geq 1-o(1)$.
Thus, these two events have a non-empty intersection and therefore it is possible to find subgraphs $H' \subset G' \subset G$ such that $G' \in \mathcal{P}$ and $\Delta(H') \leq (1/2 - \epsilon/2)np'$. Then $G' - H'$ (and hence $G - H$) must contain cycles of length $t$ for all $l \leq t \leq n$.
\end{pf}

\begin{rem}
Note that there is nothing special about the property of
``containing cycles'' and in fact if for some $\log n/n\ll p' \leq
p$ and $\alpha > 0$, we have a monotone increasing graph property
$\mathcal{Q}$ such that $G(n,p')$ a.a.s. has local resilience at
least $(\alpha+o(1))np'$ with respect to having property
$\mathcal{Q}$ then $G(n,p)$ a.a.s. has local resilience at least
$(\alpha+o(1))np$ with respect to having property $\mathcal{Q}$.
\end{rem}

From now on we may assume that $p = Cn^{-1 + 1/(l-1)}$ instead of $p \geq C n^{-1 + 1/(l-1)}$ since if we can prove the 
theorem under this condition then we can extend it to the whole range using the previous proposition. Moreover we will assume that the constant $C$ is large enough without further mentioning.

In the next two lemmas we establish some expansion properties of
random graphs.

\begin{lemma} \label{pancyclic_random_lemma10}
Fix a positive integer $l$ and $0< \epsilon <1$ and let $G=(V,E)$ be a random graph $G(n,p)$ with $p = C n^{-1+1/(l-1)}$. Then a.a.s. every subset $X \subset V$ of size $|X| \leq (2C)^{l-1}n^{-1}p^{-2}$ satisfies $(1 - \epsilon)|X|np \leq |N(X)| \leq (1 + \epsilon)|X|np$.
\end{lemma}
\begin{pf}
Fix a set $X \subset V$ of size $|X| \leq (2C)^{l-1}n^{-1}p^{-2}$.
For each $v \in V$ let $Y_v$ be indicator random variable of the
event that $v \in N(X)$. Since $|X|p=o(1)$, we have $P(Y_v = 1) = 1
- (1-p)^{|X|} = (1+o(1))|X|p$. Consider the random variable
$Y=\sum_{v \in V \backslash X} Y_v= |N(X)\setminus X|$ and note that
\[ \mathbb{E}[Y] = \sum_{v \in V \backslash X} P(Y_v =1) = (n-|X|)(1+o(1))|X|p = (1+o(1))|X|np. \]
 Moreover if $v \in V \backslash X$ then $Y_v$ are mutually independent so we can apply the Chernoff inequality to get
\[ P\big( |Y - \mathbb{E}[Y]| \geq (\epsilon/3) \mathbb{E}[Y]\big) \leq e^{-\Omega_\epsilon(\mathbb{E}[Y])}. \]
Combine this with the estimate on $\mathbb{E}[Y]$ and we have,
\[ P\big( (1 - 2\epsilon/3) |X|np \leq Y \leq (1 + 2\epsilon/3) |X|np \big) \geq 1-e^{-\Omega_\epsilon(\mathbb{E}[Y])} = 1-e^{-\Omega_\epsilon(|X|np)} \]
for large enough $n$. Finally, note that $|N(X) - Y| \leq |X| =
o(|X|np)$ and thus for large enough $n$, we have $(1 - \epsilon)
|X|np \leq |N(X)| \leq (1 + \epsilon)|X|np$ with probability at
least $1-e^{-\Omega_\epsilon(|X|np)}$.

Taking the union bound over all choices of $X$, we get
\begin{eqnarray}
 \sum_{1 \leq |X| \leq (2C)^{l-1}n^{-1}p^{-2}} e^{-\Omega_\epsilon(|X|np)} &=& \sum_{1 \leq k \leq (2C)^{l-1}n^{-1}p^{-2}} \binom{n}{k} e^{-\Omega_\epsilon(knp)}
 \leq \sum_{1 \leq k \leq (2C)^{l-1}n^{-1}p^{-2}} \left(\frac{en}{k} e^{-\Omega_\epsilon(np)}\right)^k \nonumber \\
 &\leq& \sum_{1 \leq k \leq (2C)^{l-1}n^{-1}p^{-2}} ne^{-\Omega_\epsilon(np)} \leq n^2e^{-\Omega_\epsilon(np)} = o(1). \nonumber
\end{eqnarray}
This implies the assertion of the lemma.
\end{pf}

\begin{lemma} \label{pancyclic_random_lemma1}
Fix a positive integer $l$ and $0< \epsilon <1$ and let $G=G(n,p)$
be a random graph with $p = C n^{-1+1/(l-1)}$. Then a.a.s. $G$ has
the following property . If $H$ is a subgraph of $G$ with maximum
degree at most $(1/2 - \epsilon)np$ and $G' = G-H$, then every set
$X$ with $|X| \geq 2^{-l}(np)^{l-2}$ satisfies $|N_{G'}(X)| \geq
(1/2 + \epsilon/2)n$.
\end{lemma}
\begin{pf}
It is enough to show that a.a.s. for any $H$ as above the claim
holds for every set $X$ with size exactly $|X| = 2^{-l}(np)^{l-2}$.
Fix a set $X$ of size $2^{-l}(np)^{l-2}$ and let $Y \subset V$ be a
set of size $|Y| \geq (1/2 - \epsilon/2)n$ disjoint from $X$. Then
we have $\mathbb{E}[e_G(X,Y)] = |X||Y|p
>2^{-l-2}(np)^{l-1}=2^{-l-2}C^{l-1}n$ and by the Chernoff
inequality,
\begin{eqnarray}
P\Big(\big|e_G(X, Y) - |X||Y|p\big| \geq (\epsilon/4)|X||Y|p\Big)&<& e^{-\Omega_\epsilon(|X||Y|p)} \leq e^{-\Omega_\epsilon(2^{-l-2}C^{l-1}n)}. \label{pancyclic_lemma1_eqn3}
\end{eqnarray}
Thus with probability at least the right hand side of
(\ref{pancyclic_lemma1_eqn3}) we have $e_G(X,Y) \geq (1 -
\epsilon/4)|X||Y|p \geq (\frac{1}{2} - 3\epsilon/4)|X|np$. Since
there are at most $2^{2n}$ possible choices of the pairs $X, Y$ and
the right hand side of (\ref{pancyclic_lemma1_eqn3}) is $\ll
2^{-2n}$ for large enough $C$, we a.a.s. have $e_G(X,Y) >
(\frac{1}{2} - \epsilon)|X|np$ for every pair $X,Y$ as above.

On the other hand, we know that $e_H(X,Y) \leq (1/2 -
\epsilon)np|X|$. Therefore a.a.s. $e_{G'}(X,Y) \geq e_G(X,Y) -
e_H(X,Y) > 0$. This implies that $N_{G'}(X) \cap Y \neq \emptyset$
for all $Y$ with $|Y| \geq (1/2 - \epsilon/2)n$. Thus $|N_{G'}(X)|
\geq n - (1/2 - \epsilon/2)n = (1/2 + \epsilon/2)n$.
\end{pf}

We also need the following lemma that proves expansion property for
subgraphs of $G(n,p)$ with large minimum degree.

\begin{lemma} \label{lemma_pancyclic2}
If $p = C n^{-1+1/(l-1)}$ and $\epsilon' > 0$ then a.a.s.
every subgraph $G' \subset G(n,p)$ with minimum degree at
least $\epsilon' np$ satisfies the following expansion property.
For all $X \subset V$ with $|X| \leq \frac{1}{80}\epsilon' n$, $|N_{G'}(X)\backslash X| \geq
2|X|$.
\end{lemma}
\begin{pf}
Assume to the contrary that there exists a set $X \subset V$ such
that $|X| \leq \frac{1}{80}\epsilon' n$ and $|N_{G'}(X) \backslash
X| < 2|X|$, and let $Y = X \cup N_{G'}(X)$ so that $|Y| \leq 3|X|
\leq \frac{1}{20} \epsilon' n$. Then by the minimum degree condition
we know that $e_{G'}(Y) \geq \frac{1}{2} |X| \epsilon' np \geq
\frac{1}{8}|Y| \epsilon' np$. Now we will estimate the probability
that such event can happen for a set $Y$ with $|Y| = a$. We can
restrict the range to $\epsilon' np \leq a \leq \frac{1}{20}
\epsilon' n$ since $G'$ has minimum degree at least $\epsilon' np$.
The probability that there exists a set of size $a$ which spans at
least $\frac{1}{8}a \epsilon' np$ edges is,
\begin{eqnarray}
 \binom{n}{a} \binom{a(a-1)/2}{a\epsilon'np/8} p^{\epsilon' a np/8 }
  &\leq& \left(\frac{en}{a}\right)^a \left(\frac{4ea}{\epsilon' n p}\right)^{\epsilon' a n p/8}  p^{\epsilon' a np/8 }
  = \bigg(\frac{en}{a} \Big(\frac{4ea}{\epsilon' n }\Big)^{\epsilon' np/8} \bigg)^{a}  \nonumber \\
  &\ll& \bigg( \Big(\frac{e}{4}\Big)^{\epsilon' np/8} \bigg)^{a} \ll n^{-2}\,. \nonumber
\end{eqnarray}
Summing over all $\epsilon' np \leq a \leq \frac{1}{20}\epsilon' n$ we get that the probability that there is a set violating the assertion of the lemma is $o(1)$.
\end{pf}

\section{Proof of Theorem \ref{pancyclic_thm4} }
\label{proofofrandomgraphthm_section}

In this section we prove Theorem \ref{pancyclic_thm4}. First we need an additional lemma which gives us more properties of a random graph with deleted edges.


\begin{lemma}  \label{lemma_pancyclic1} For every integer $l\geq 3$ and $\epsilon>0$ there exists $C=C(\epsilon)$ such that if $p = C n^{-1+1/(l-1)}$ then $G=G(n,p)$ a.a.s. has the following properties. Let $H$ be a subgraph of $G$ with maximum degree at most $(\frac{1}{2} - \epsilon)np$, $G'=G-H$ and $v \in V$, then
\begin{itemize}
\item[(a)] For every $1 \leq i \leq l-2$, $2^{-i} (np)^{i} \leq |N^{(i)}_{G'}(v)| \leq (1 + \epsilon)^i(np)^{i}$.
\item[(b)] $|N^{(l-1)}_{G'}(v)| \geq (\frac{1}{2} + \frac{\epsilon}{3})n$.
\item[(c)] For every vertex $w \in V$ whose distance from $v$ is at least $l-2$, $|N^{(l-2)}_G(v) \cap N_G(w)| \leq \log n$.
\end{itemize}
\end{lemma}
\begin{pf}
Let $v$ be a vertex of $G$ and for simplicity of notation let $Y_j = N_{G'}^{(j)}(v)$ for $j=1,\ldots, l-2$.

(a) By using induction we will show that $2^{-i} (np)^{i} \leq |Y_i|
\leq (1+\epsilon)^i(np)^{i}$ for all $1 \leq i \leq l-2$. For the
initial case $i=1$ by Lemma \ref{pancyclic_random_lemma10} we have
$|Y_1| \leq |N^{(1)}_{G}(v)| \leq (1 + \epsilon)np$. By the same
lemma we also know that $|N_G^{(1)}(v)| \geq (1 - \epsilon)np$.
Therefore $|Y_1| \geq |N^{(1)}_{G}(v)| - |N^{(1)}_{H}(v)| \geq
np/2$. Now assume that we have established the claim up to some $i
\leq l-3$ and let us look at the case $i+1$. First notice
$(1+\epsilon)^{l-3}(np)^{l-3} \leq (2C)^{l-1}n^{-1}p^{-2}$ so that
we can apply Lemma \ref{pancyclic_random_lemma10}. Then the upper
bound easily follows as $|Y_{i+1}| \leq |N_{G}(Y_i)| \leq
(1+\epsilon)np|Y_i| \leq (1+\epsilon)^{i+1}(np)^{i+1}$ by that lemma
and the inductive hypothesis. To obtain the lower bound, we use that
$|N_{H}(Y_i)| \leq \Delta(H)|Y_i|$ and that $|N_{G}(Y_i)| \geq (1 -
\epsilon/2)np|Y_i|$ by Lemma \ref{pancyclic_random_lemma10} (where
we substitute $\epsilon/2$ instead of $\epsilon$). Therefore,
\[   |N_{G'}(Y_i)| \geq |N_{G}(Y_i)| - |N_{H}(Y_i)| \geq (1 - \epsilon/2)np|Y_i| - (1/2 - \epsilon)np|Y_i|
        = (1/2 + \epsilon/2)np|Y_i|.  \]
Recall the recursive formula $Y_{i+1} = N_{G'}(Y_i) - Y_i - Y_{i-1}$. By the inductive hypothesis it is easy to check that $|Y_{i-1}| = o(|Y_i|)$. Thus,
\[ |Y_{i+1}| \geq |N_{G'}(Y_i)| - |Y_i| - |Y_{i-1}| \geq  (1/2 + \epsilon/2)np|Y_i| - o(np|Y_i|) \geq 2^{-i-1} (np)^{i+1}\,, \]
 which completes the proof of the first part.

(b) By part (a) we have $2^{-l+2}(np)^{l-2} \leq |Y_{l-2}| \leq (1+\epsilon)^{l-2}(np)^{l-2}$ and $|Y_{l-3}| \leq (1+\epsilon)^{l-3}(np)^{l-3}$. Apply Lemma \ref{pancyclic_random_lemma1} to get $|N_{G'}(Y_{l-2})| \geq (1/2 + \epsilon/2)n$. Then
$$|Y_{l-1}| \geq |N_{G'}(Y_{l-2})| - |Y_{l-2}| - |Y_{l-3}| \geq (1/2 + \epsilon/2)n - (2np)^{l-2}\,,$$ and therefore for large enough $n$, $|Y_{l-1}| \geq (1/2 + \epsilon/3)n$.

(c) Condition on the event that $|N^{(i)}_{G}(v)| \leq
(1+\epsilon)^{i}(np)^{i}$ for all $i=0,\ldots,l-2$ and let $X =
\cup_{i=0}^{l-3} N^{(i)}_{G}(v)$.  Notice that so far we only
exposed the edges inside $G[X]$ and the edges connecting $X$ to
$N^{(l-2)}_{G}(v)$. Therefore for any vertex $w \notin X$ which is
at distance is at least $l-2$ from $v$, the edges between $w$ and
$N_{G}^{(l-2)}(v)$ are not yet exposed. Thus we can bound the
probability that $w$ has degree at least $\log n$ in
$N^{(l-2)}_{G}(v)$ as follows:
\begin{eqnarray}
\binom{|N^{(l-2)}_{G}(v)|}{\log n} p^{\log n} &<& \bigg(\frac{e|N^{(l-2)}_{G}(v)|p}{\log n}\bigg)^{\log n}
< \bigg(\frac{e2^{l-2}(np)^{l-2}p}{\log n}\bigg)^{\log n}
  = \bigg(\frac{e2^{l-2}C^{l-1}}{\log n}\bigg)^{\log n}. \nonumber
\end{eqnarray}
Since the last estimate is $o(n^{-2})$, a.a.s. every pair of
vertices as above satisfies the claim.
\end{pf}

Now we are ready to prove Theorem \ref{pancyclic_thm4}. First we
restate it in a more accurate and general form.

\begin{thm}
For every $\epsilon>0$ there exists $C$ such that if $p \geq Cn^{-1+1/(l-1)}$ then $G(n,p)$ almost surely has local resilience at least $(1/2 - \epsilon)np$ with respect to being pancyclic.
\end{thm}
\begin{pf}
By Proposition \ref{probrestriction_prop} we may assume that $p = Cn^{-1+1/(l-1)}$
where $C$ is taken to be the maximal of the corresponding constants in Theorem \ref{pancyclic_thm2} and Lemma \ref{lemma_pancyclic1}. Let $G=G(n,p)$, $H$ be a subgraph of maximum degree $\Delta(H) \leq (\frac{1}{2} - \epsilon)np$, and $G' = G - H$. The proof consists of three parts. In each part we will show the existence of short, medium length, and long cycles, respectively, in $G'$.\\

\noindent \textbf{Short Cycles.}\,
The existence of cycles of length $l$ to $2l-2$ in $G'$ is a direct corollary of Haxell, Kohayakawa and {\L}uczak's Theorem \ref{pancyclic_thm2}.\\

\noindent \textbf{Medium Length Cycles.}\, Now we show the existence
of cycles of length $2l-1$ up to $\frac{1}{320}\epsilon n$. Fix a
vertex $v \in V$ and let $Y=N_{G'}^{(l-1)}(v)$. Then by Lemma
\ref{lemma_pancyclic1} part (b) a.a.s.  $|Y| \geq (\frac{1}{2} +
\epsilon/3)n$. By applying the Chernoff inequality and then taking
the union bound over sets $Y$ of appropriate sizes, we know that
a.a.s. $e_G(Y) \geq (1 - \epsilon/6)\binom{|Y|}{2}p \geq
\frac{1}{4}|Y|np$. And by the restriction on the maximum degree of
$H$ we know that $e_H(Y) \leq \frac{1}{2}(\frac{1}{2} -
\epsilon)|Y|np$. Therefore,
\[ e_{G'}(Y) \geq e_{G}(Y) - e_{H}(Y) \geq \frac{1}{2}\epsilon|Y|np \geq \frac{1}{4}\epsilon n^2p. \]
Thus by Lemma \ref{lemma_findingmindegree}, we can find a subgraph
$G_1 \subset G'[Y]$ with minimum degree at least
$\frac{1}{4}\epsilon np$. Fix any vertex $v_{l-1} \in V(G_1)$ and
let $v_i \in N^{(i)}_{G'}(v)$ for $i= 1,\ldots, l-2$ be the vertices
of a path $vv_1 v_2 \ldots v_{l-1}$ in $G'$ from $v$ to $v_{l-1}$.
Delete every vertex in $N^{(l-2)}_{G'}(v_1) \cap N^{(l-1)}_{G'}(v)
\subset N^{(l-1)}_{G'}(v)$ except $v_{l-1}$ from $G_1$ to obtain
$G_2$. Then by Lemma \ref{lemma_pancyclic1} part (c), $\delta(G_2)
\geq \delta(G_1) - \log n$, and so for large enough $n$, $G_2$ has
minimum degree at least $\frac{1}{8}\epsilon np$. Now by Lemma
\ref{lemma_pancyclic2}, $G_2$ has the property that every subset $X$
of size $|X| \leq \frac{1}{640}\epsilon n$ satisfies
$|N_{G_2}(X)\backslash X| \geq 2|X|$. Therefore by P\`{o}sa's
rotation-extension Lemma \ref{lemma_pancyclic3} we can find a path
$P$ of length at least $\frac{1}{320}\epsilon n$ starting at
$v_{l-1}$ inside $G_2$. Let this path be $P=v_{l-1} w_1 \ldots
w_{\epsilon'n}$ where $\epsilon' \geq \frac{1}{320}\epsilon$.
Finally observe that for any vertex $w_t (t>0)$ there is a path
$vz_1z_2 \ldots z_{l-2}w_t$ in $G'$ such that $z_j \in
N^{(j)}_{G'}(v)$ for $j=1, \ldots l-2$.
Moreover since we deleted vertices that can be reached from $v_1$, $v_j \neq z_j$ for all $1 \leq j \leq l-2$. Thus we have a cycle $vv_1\ldots v_{l-1}w_1\ldots w_tz_{l-2}\ldots z_1 v$ which has length $t + 2l - 2$. Since $t$ can be arbitrarily chosen in the range $1 \leq t \leq \epsilon' n$, we are done with the second part of the proof.\\

\noindent \textbf{Long Cycles.}\, Let $\alpha =
\frac{1}{320}\epsilon$. In this part we will show how to find all
cycles of length from $\alpha n$ to $n$ in $G'$. For a fixed integer
$n^*$ satisfying $\alpha n \leq n^* \leq n$ choose uniformly at
random $n^*$ vertices $V^{*}$ out of $V$ and let $G^{*}=G[V^{*}],
H^{*}=H[V^{*}]$. Let $\mathcal{P}$ be the graph property of a graph
on $n^*$ vertices having local resilience at least
$(1/2-\epsilon/2)n^*p$ with respect to hamiltonicity. We claim that
with probability $1-o(n^{-1})$, $P(G^{*} \in \mathcal{P} |G) \geq
1/2$. First note that $G^{*}$ has distribution $G(n^*, p)$ and apply
Sudakov and Vu's Theorem \ref{pancyclic_thm1} to get $P(G^{*} \notin
\mathcal{P}) = o(n^{-1})$. Let $A$ be the event in the probability
space $G(n,p)$ such that $P(G^{*} \notin \mathcal{P} |G) \geq 1/2$.
Then we have $(1/2)P(A) \leq P(G^{*} \notin \mathcal{P}) =
o(n^{-1})$. Therefore $P(A) = o(n^{-1})$, or in other words $P(G^{*}
\in \mathcal{P} |G) \geq 1/2$ with probability at least
$1-o(n^{-1})$. Let $\mathcal{A}_{n^{*}}$ be the collection of graphs
$G$ having this property.

On the other hand, observe that the degree of a vertex in $H^{*}$
follows the hypergeometric distribution and thus we can apply Lemma
\ref{lemma_hypergeometric}. Hence for a vertex $v \in V^{*}$,
\[ P\big(|\deg_{H^{*}}(v) - (1/2 - \epsilon)n^*p| \geq \epsilon n^*p /2 \big) \leq e^{ -\Omega_{\epsilon}(n^*p) } \leq e^{-\Omega_{\epsilon} (np)}\,, \]
thus a.a.s. every vertex in $V^{*}$ has degree at most $(1/2 -
\epsilon/2)n^* p$ in $H^{*}$. We can conclude that if $G \in
\mathcal{A}_{n^*}$ then there exists a set $V^{*}$ of size $n^*$
such that $G^{*} \in \mathcal{P}$ and $\Delta(H^{*}) \leq (1/2 -
\epsilon/2)n^*p$. This gives a hamilton cycle inside $G^{*} - H^{*}$
which is a cycle of length $n^*$ inside $G' = G-H$.

Finally note that since $G \in \mathcal{A}_{n^*}$ with probability
at least $1 - o(n^{-1})$, cycles of length $n^*$ exist with
probability at least $1-o(n^{-1})$ for any fixed $n^*$ by the
previous observation. Therefore by taking the union bound we can see
that a.a.s. $G'$ simultaneously contains cycles of length $n^*$ for
all $\alpha n \leq n^* \leq n$. This concludes the proof.
\end{pf}


\section{Properties of pseudo-random graphs}
\label{pseudorandomgraphproperty_section}

Here we collect properties of pseudo-random graphs which we will use
later to prove Theorem \ref{pancyclic_pseudo_thm01}. The main fact
that we use about $(n,d,\lambda)$-graphs is the following formula
established by N. Alon (see, e.g., \cite{MR2223394}) which connects
between eigenvalues and edge distribution.
\begin{lemma} \label{pancyclic_pseudo_lemma0}
If $G=(V,E)$ is an $(n,d,\lambda)$-graph, then for any $X,Y \subset V$ we have,
\[ \Big|e(X,Y) - \frac{d}{n}|X||Y|\Big| \leq \lambda \sqrt{|X||Y|}. \]
\end{lemma}

As in Section \ref{randomgraphproperty_section} we will prove
several lemmas that establish some expansion properties of
pseudo-random graphs. These lemmas correspond to Lemma
\ref{pancyclic_random_lemma10}, Lemma \ref{pancyclic_random_lemma1}
and Lemma \ref{lemma_pancyclic2} in the random graph case.

\begin{lemma} \label{pancyclic_pseudo_lemma10}
Let $\epsilon, \epsilon'$ be such that $0 \leq 5\epsilon' <
\epsilon, k \geq 3$, and let $G = (V,E)$ be an
$(n,\epsilon',d,\lambda)$-graph with $d^{k-1} /n = \omega(n)
\lambda^{k-2}$ where $\omega(n) \rightarrow \infty$. Then $G$ has
the following property. If $H$ is a subgraph of $G$ with $\Delta (H)
\leq (1/2 - \epsilon) d$ and $G' = G-H$, then every set $X$ with
$|X| \leq \epsilon n/4$ satisfies $|N_{G'}(X)| \geq \min
\big(\epsilon n /2, \frac{d^2}{4\lambda^2} |X|\big)$.
\end{lemma}
\begin{pf}
Let $Y = N_{G'}(X)$ and assume that $|Y| \leq \epsilon n/2$ as otherwise we are done. Since $G'$ has minimum degree at least $(1-\epsilon')d - \Delta(H) \geq (1-\epsilon')d - (1/2-\epsilon)d \geq (1/2 + 4\epsilon/5)d$, we have
 \begin{eqnarray}
  e_{G'}(X, Y) &\geq& (1/2 + 4\epsilon/5)d|X| - 2e_G(X)
                   \geq (1/2 + 4\epsilon/5)d|X| - (d|X|^2/n + \lambda |X|) \nonumber \\
                   &\geq& (1/2 + 4\epsilon/5)d|X| - (\epsilon d/4 + \lambda) |X|                   \geq (1/2 + \epsilon/2)d|X|. \label{eq:111}
 \end{eqnarray}
On the other hand, since $|Y| \leq \epsilon n/2$, we have:
 \begin{eqnarray}
  e_G(X, Y) &\leq& \frac{d|X| |Y|}{n} + \lambda \sqrt{|X||Y|}
                   \leq (\epsilon/2) d |X| + \lambda \sqrt{|X||Y|}. \label{eq:112}
 \end{eqnarray}
Therefore by (\ref{eq:111}),(\ref{eq:112}) and $e_{G'}(X,Y) \leq e_{G}(X,Y)$ we have, $(1/2 + \epsilon/2)d|X| \leq (\epsilon/2) d |X| + \lambda \sqrt{|X||Y|}$ which implies $|Y| \geq \frac{d^2}{4\lambda^2} |X|$.
\end{pf}

\begin{lemma} \label{pancyclic_pseudo_lemma1}
Let $k \geq 3$ and let $G = (V,E)$ be an $(n,\epsilon',d,\lambda)$-graph with $d^{k-1} /n = \omega(n) \lambda^{k-2}$ where $\omega(n) \rightarrow \infty$. Then for any function $\delta=\delta(n)$ such that $1 \ll \delta \ll d^2/\lambda^2$, $G$ has the following property. If $H$ is a subgraph of $G$ with $\Delta (H) \leq (1/2 - \epsilon) d$ and $G' = G-H$, then every set $X$ with $|X| \geq \delta (\lambda^2 / d^2)n$ satisfies $|N_{G'}(X)| \geq (1/2 + \epsilon/2)n$.
\end{lemma}
\begin{pf}
We only have to verify this for sets of size exactly
$\delta(\lambda^2 / d^2)n$, so assume for the contrary that there
exists $X \subset V$ of size $|X| = \delta(\lambda^2 / d^2)n$ which
has $|N_{G'}(X)| < (1/2 + \epsilon/2)n$ and define $Y = V \backslash
(X \cup N_{G'}(X))$. Since $|X|=o(n)$, we have $|Y| > (1/2 -
2\epsilon/3)n$. Therefore by Lemma \ref{pancyclic_pseudo_lemma0},
\begin{eqnarray}
 e_G(X, Y) \geq \frac{d|X||Y|}{n} - \lambda \sqrt{|X||Y|} \geq (1/2 - 2\epsilon/3)\delta (\lambda^2/d) n - \delta(\lambda^2/d) n > (1/2 - \epsilon) \delta (\lambda^2/d) n. \nonumber
\end{eqnarray}
On the other hand by the maximum degree restriction, $e_H(X,Y) \leq (1/2 - \epsilon)d|X| = (1/2 - \epsilon)\delta (\lambda^2 /d) n$. But since there are no edges between $X$ and $Y$ we must have $0=e_{G'}(X,Y) \geq e_G(X,Y) - e_H(X,Y) >  0$ which gives us a contradiction.
\end{pf}

The next lemma proves expansion property for subgraphs of
$(n,d,\lambda)$-graphs with large minimum degree.

\begin{lemma} \label{pancyclic_pseudo_lemma2}
Let $G=(V,E)$ be an $(n,d,\lambda)$-graph with $\lambda =o(d)$, and
let $G'$ be a subgraph of $G$ with  $\delta(G') \geq \epsilon d$ for
some fixed constant $\epsilon > 0$. Then every $X \subset V(G')$
with $|X| \leq \epsilon n /10$ satisfies $|N(X) \backslash X| \geq
2|X|$.
\end{lemma}
\begin{pf}
Assume to the contrary that there exists a set $X \subset V(G')$
with $|X| \leq \epsilon n/10$ and $|N(X) \backslash X| < 2|X|$. Let
$A = X \cup N(X)$ and note that $|A|< 3|X|$. Then by Lemma
\ref{pancyclic_pseudo_lemma0},
\[ e_G(A)=e_G(A,A)/2 \leq \frac{d}{2n}|A|^2 + \frac{\lambda}{2}|A| \leq |X| (9\epsilon d/10 + 3\lambda)/2.  \]
On the other hand, since $G'$ has minimum degree at least $\epsilon
d$, we have
\[ e_G(A) \geq e_{G'}(A) \geq |X|\epsilon d/2 \]
which is a contradiction, since $\lambda = o(d)$.
\end{pf}


\section{Proof of Theorem \ref{pancyclic_pseudo_thm01}}
\label{proofofpseudorandomgraphthm_section}

In this section we prove Theorem \ref{pancyclic_pseudo_thm01}. As in the random graph case, we need an additional lemma which gives us more properties of a pseudo-random graph with deleted edges.

\begin{lemma} \label{pancyclic_pseudo_lemma11}
Fix $\epsilon > 0, k \geq 3$ and let $G = (V,E)$ be a
$(n,d,\lambda)$-graph with $d^{k-1} /n = \omega(n) \lambda^{k-2}$
where $\omega(n) \rightarrow \infty$. Let $H$ be a subgraph of $G$
with $\Delta (H) \leq (1/2 - \epsilon) d$, $G'=G-H$ and $v \in V$.
Then there exist $l, 1 \leq l \leq \lfloor(k-1)/2\rfloor$ and sets
$X_i(v), Y_i(v)$ for  $i=0, 1, \ldots, l$ such that,
\begin{itemize}
\item[(a)] $X_0(v) = Y_0(v) = \{ v \}$, $X_i(v) \cap Y_i(v) = \emptyset, |X_i(v)| = |Y_i(v)|$ for all $i \neq
0$;
\item[(b)] $|X_{i+1}(v)|\geq \frac{d^2}{16\lambda^2}|X_i(v)|$, $|Y_{i+1}(v)|\geq \frac{d^2}{16\lambda^2}|Y_i(v)|$ for all $i=0, \ldots l-2$
 and $|X_{i}(v)|=|Y_{i}(v)| \leq (\epsilon/(8k))n$ for all $i=0,1 \ldots,
 l$;
\item[(c)] Let $Z_i(v) = \cup_{j=0}^{i} (X_j(v) \cup Y_j(v)) $. Then $X_{i+1}(v) \subset N_{G'}(X_i(v)) \backslash Z_i(v)$, $Y_{i+1}(v) \subset N_{G'}(Y_i(v)) \backslash Z_i(v)$ for all $0 \leq i \leq l-1$.
\item[(d)] $|X_{l}(v)|=|Y_{l}(v)| \geq \delta (\lambda^2 / d^2) n$ for some function $\delta=\delta(n) \rightarrow \infty$.
\end{itemize}
\end{lemma}
\begin{pf}
Let $\delta = \delta(n)=\min(d / \lambda, (\omega(n))^{1/2})$ and
note that indeed $\delta \rightarrow \infty$. Given a vertex $v \in
V$, we will inductively construct sets $X_i = X_i(v), Y_i = Y_i(v)$
satisfying the condition above. Since $G'$ has minimum degree at
least $(1/2 + \epsilon)d$, put $d/4$ vertices of $N_{G'}(v)$ into
$X_1$ and put another $d/4$ vertices into $Y_1$. Suppose that for
some $i \leq \lfloor (k-1)/2 \rfloor - 1$ we have already
constructed $X_0, Y_0, \ldots, X_i, Y_i$ satisfying conditions
$(a),(b),(c),(d)$. Next we show how to construct $X_{i+1}, Y_{i+1}$.
Let $Z_i = \cup_{j=0}^{i} (X_j \cup Y_j)$ and note that $\big| Z_i
\big| \leq 2i|X_i| \leq k|X_i|$. If $|X_i| = |Y_i| \geq \delta
(\lambda^2 / d^2) n$ then define $l=i$ and stop the process.
Otherwise $|X_i| = |Y_i| < \delta (\lambda^2 / d^2) n \leq (\lambda
/d)n = o(n)$ and by Lemma \ref{pancyclic_pseudo_lemma10} we have
that $|N(X_i)|, |N(Y_i)| \geq \min (\epsilon n /2,
\frac{d^2}{4\lambda^2} |X|)$. Thus
\[ \big| N(X_i) \backslash Z_i \big| \geq 
\min\Big(\epsilon n /2 -  k|X_i|, \frac{d^2}{4\lambda^2} |X_i| - k|X_i| \Big) \geq \min \Big(\epsilon n /4 , 
\frac{d^2}{8\lambda^2}|X_i|\Big) \]
and a similar inequality also holds for $N(Y_i)$. Therefore, by
splitting the vertices of $N(X_i) \cap N(Y_i)$ between $X_{i+1}$ and
$Y_{i+1}$ we can always choose $X_{i+1} \subset N(X_i) \backslash
Z_i$ and $Y_{i+1} \subset N(Y_i) \backslash Z_i$ so that $X_{i+1}
\cap Y_{i+1} = \emptyset$ and $|X_{i+1}|=|Y_{i+1}| \geq (1/2) \min
(\epsilon n /4 , \frac{d^2}{8\lambda^2}|X_i|)$. If $\epsilon n /4
\leq \frac{d^2}{8\lambda^2} |X_i|$ then $|X_{i+1}|, |Y_{i+1}| \geq
\epsilon n / 8$, so stop the process and define $l=i+1$. Otherwise
we can make $|X_{i+1}|, |Y_{i+1}| \geq \frac{d^2}{16\lambda^2}
|X_i|$ and continue. Note that (b) holds in this case. If the
process does not terminate after constructing $X_1, \ldots,
X_{\lfloor (k-1)/2 \rfloor}$ and $Y_1, \ldots, Y_{\lfloor (k-1)/2
\rfloor}$ then by property (b) we get that $\frac{d}{4}
(\frac{d^2}{16\lambda^2})^{\lfloor (k-1)/2 \rfloor - 1} \leq
X_{\lfloor (k-1)/2 \rfloor} < \delta \frac{\lambda^2}{d^2} n$. This
implies:

\[ \delta \frac{\lambda^2}{d^2} n > \frac{d}{4}\left(\frac{d^2}{16\lambda^2}\right)^{\lfloor (k-1)/2 \rfloor - 1} \geq \frac{d}{4}\left(\frac{d^2}{16\lambda^2}\right)^{k/2 - 2} = \frac{d^{k-3}}{4^{k-3}\lambda^{k-4}} = \frac{\omega(n)}{4^{k-5}} \cdot \frac{\lambda^2}{d^2} n\]
which is a contradiction, since $\delta \ll \omega(n)$. Finally note
that we can always shrink final sets $X_l, Y_l$ so that they become
smaller than $(\epsilon /(8k))n$. Since  $|X_{l-1}|=|Y_{l-1}| <
\delta (\lambda^2/d^2)n \ll (\epsilon/(8k)) n$, (b) holds for all
$i=0,1 \ldots, l$. Thus we can find sets as claimed.
\end{pf}

We are ready to prove Theorem \ref{pancyclic_pseudo_thm01}. First we restate it here with more quantifiers.

\begin{thm}
Fix $\epsilon>0$ and let $k$ be either 3 or an even integer
satisfying $k \geq 4$, and let $G=(V,E)$ be a $(n,d,\lambda)$-graph
satisfying $d^{k-1} /n = \omega(n) \lambda^{k-2}$ where $\omega(n)
\rightarrow \infty$. Then for large enough $n$, $G$ has local
resilience at least $(1/2 - \epsilon)d$ with respect to containing
cycles of length $t$ for $k \leq t \leq n$.
\end{thm}
\begin{pf}
Let $H$ be a subgraph of $G$ with $\Delta(H) \leq (1/2 - \epsilon)d$
and let $G' = G-H$. If $d>(1-\epsilon)n$, then $G'$ has minimum
degree larger than $(1/2+\epsilon)d>(1/2+\epsilon)(1-\epsilon)n>n/2$
for small enough $\epsilon>0$. Hence by Bondy's theorem, mentioned
in the introduction, $G'$ is pancyclic. Assume therefore that
$d\leq(1-\epsilon)n$. This implies that $\lambda=\Omega(\sqrt{d})$.
Indeed, let $A$ be the adjacency matrix of $G$ and $d=\lambda_1,
\ldots, \lambda_n$ be its eigenvalues. The trace of $A^2$ is the
number of ones in $A$, which implies that
$$ 2|E(G)|=nd=Tr(A^2)=\sum_{i=1}^n\lambda_i^2\le d^2+(n-1)\lambda^2.$$
Solving the above inequality for $\lambda$ establishes the claim.

A proof of the theorem consists of three parts. In each part we will show the existence of short, medium length and long cycles respectively.\\

\noindent \textbf{Short Cycles.}\, For $k=3$, the existence of
$3$-cycles is a direct corollary of Sudakov, Szab\'{o} and  Vu's
Theorem \ref{pancyclic_pseudo_thm1}. Also in this case we have that
$d^3/\lambda^2\geq d^2/\lambda \gg n$. Therefore for $k=3$ the
existence of cycles of length $4, \ldots, n$ follows from the proof
of case $k=4$. So from now on we assume that $k=2k'$ is even.
Since $\lambda = \Omega(\sqrt{d})$, we have $d^{k-1} = \omega(n) \lambda^{k-2} n = \omega(n) \Omega(d^{k/2 - 1}) n$. Therefore $d \gg  n^{2/k}$ and by Bondy and Simonovits' Theorem \ref{bondysimonovits1_thm} $G'$ must have a $k$-cycle.\\

\noindent \textbf{Medium Length Cycles.}\, The next step is to prove
the existence of cycles of length  from $k+1$ up to $\epsilon n/
20$.  Fix a vertex $v$ and apply Lemma
\ref{pancyclic_pseudo_lemma11} to find sets $X_1=X_1(v), \ldots,
X_l=X_l(v), Y_1=Y_1(v), \ldots, Y_l=Y_l(v)$ where $l \leq k'-1$ with
$|X_{l}|=|Y_{l}| \geq \delta (\lambda^2 / d^2) n$ and
$|X_{i}|=|Y_{i}| \leq (\epsilon/(8k))n$ for all $i=1, \ldots, l$.
Then $|\cup_{i=0}^{l}X_{i} \cup Y_{i}| \leq \epsilon n /4$. By Lemma
\ref{pancyclic_pseudo_lemma1} we know that $|N_{G'}(X_l)| \geq (1/2
+ \epsilon/2)n$ and so if we let $Z = N_{G'}(X_l) \backslash
\big(\cup_{i=0}^{l}X_{i} \cup Y_{i}\big)$ then $|Z| \geq (1/2 +
\epsilon/4)n$. Since $\lambda=o(d)$, by Lemma
\ref{pancyclic_pseudo_lemma0}, we have $e_G(Z) \geq d|Z|^2/(2n) -
\lambda|Z|/2 \geq d|Z|/4$. On the other hand $e_H(Z) \leq
(1/2-\epsilon)d|Z|/2$. Hence  $e_{G'}(Z) \geq e_{G}(Z) - e_{H}(Z)
\geq \epsilon d|Z|/2$. This implies by Lemma
\ref{lemma_findingmindegree} that inside $G'[Z]$ we can find a
subgraph $G_1 \subset G'$ which has minimum degree at least
$\epsilon d/2$. Then using Lemma \ref{pancyclic_pseudo_lemma0} it is
easy to check that $|V(G_1)| \geq \epsilon n /3$ and so we can
choose a set $W \subset V(G_1) \subset Z$ of size $\epsilon n / 8$.
Then by Lemma \ref{pancyclic_pseudo_lemma1}, we have that both
$|N_{G'}(W)|, |N_{G'}(Y_l)| \geq (1/2 + \epsilon)n$. Therefore the
set $\big(N_{G'}(W) \cap N_{G'}(Y_l)\big) \backslash
\big(\cup_{i=0}^{l}X_{i} \cup Y_{i}\big)$ has size at least
$2\epsilon n-\epsilon n/4>0$. In particular, there must exist a
vertex $p \in \big(N_{G'}(W) \cap N_{G'}(Y_l)\big) \backslash
\big(\cup_{i=0}^{l}X_{i} \cup Y_{i}\big)$, and let $y_l \in Y_l$, $w
\in W$ be neighbors of $p$ in $G'$.
Since $y_l \in Y_l$, by the definition of $Y_l$, there exists a path $vy_1y_2 \ldots y_l$ in $G'$ from $v$ to $y_l$ such that $y_i \in Y_i$ for $i=1, \ldots l$. If $p \in V(G_1)$ then let $G_2 = G_1 \backslash \{ p\}$ and note that $G_2$ has minimum degree at least $\epsilon d / 4$. Now by Lemma \ref{pancyclic_pseudo_lemma2} every set $X \subset V(G_2)$ of size $|X| \leq \epsilon d /40$ satisfies $|N_{G_2}(X)\backslash X| \geq 2|X|$. Then by P\`{o}sa's rotation-extension Lemma \ref{lemma_pancyclic3} we know that there exists a path $P =  v_0 v_1 \ldots v_t $ starting at $v_0 = w$ which  has length at least $t \geq \epsilon n/ 20$ inside $G_2$. For an arbitrary $v_i \in P, i \geq 0$ since $v_i \in V(G_2) \subset N_{G'}(X_l)$, there is a path $vx_1\ldots x_lv_i$ in $G'$ such that $x_i \in X_i$ for $i=1,\ldots,l$. Thus we have a cycle $v x_1 x_2 \ldots x_l v_i v_{i-1} \ldots v_0 p y_l y_{l-1} \ldots y_1v$ which has length $2(l + 1) + i + 1 \leq k+i+1$. Since $i$ can be arbitrarily chosen in the range $0 \leq i \leq \epsilon n/20$, we are done with the second part of the proof.\\

\noindent \textbf{Long Cycles.}\, The final step is to prove the
existence of cycles of length $\epsilon n/20$ to $n$. Pick $\epsilon
/20 \leq \alpha \leq 1$ such that $n^{*} = \alpha n$ is an integer.
Let $V^{*} \subset V$ be a set of size $n^{*}$ chosen uniformly at
random and $G^{*}=G[V^{*}], H^{*}=H[V^{*}]$ be the induced subgraph
of $G, H$ respectively. Since $d \gg n^{2/k}$, by the concentration
of the hypergeometric distribution (Lemma
\ref{lemma_hypergeometric}), for every vertex $v \in V^{*}$ we have
that
\begin{eqnarray}
P\big(deg_{H^{*}}(v) \geq (1/2 - \epsilon /2)\alpha d \big) &\leq& e^{-\Omega_\epsilon(\alpha d)} \leq e^{-\Omega_\epsilon(n^{2/k})}. \nonumber
\end{eqnarray}
Similarly, with probability $1-o(n^{-1})$, the graph $G^{*}$ has
minimum degree $(1-\epsilon/6)\alpha d$ and therefore if $n$ is
large enough, for every $\epsilon n/20\le n^*\le n$ there exists a
choice of $V^{*}$ where $\Delta(H^{*}) \leq (1/2 - \epsilon/2)\alpha
d$ and $\delta(G^{*}) \geq (1-\epsilon/6)\alpha d$. Moreover since
$G^{*}$ is an induced subgraph of $G$, its edge distribution is
still governed by the estimate from Lemma
\ref{pancyclic_pseudo_lemma0}. Therefore $G^{*}$ is an $(\alpha n,
\epsilon/6, \alpha d, \lambda)$-graph. Now by Sudakov and Vu's
Theorem \ref{pancyclic_pseudo_thm2}, for large enough $n$, $G^{*}$
has local resilience at least $(1/2 - \epsilon/2)\alpha d$ with
respect to being hamiltonian. Thus for the choice of $V^{*}$ as
above, $G^{*}-H^{*}$ must contain a hamilton cycle which is a cycle
of length $n^{*}$ in $G'$. This concludes the proof.
\end{pf}


\section{Concluding Remarks}
\label{concludingremarks_section}

\subsection{Cycles through a given vertex}

In this paper we found cycles inside a subgraph of random and
pseudo-random graphs. We proved that under certain conditions there
exist cycles of various lengths \textit{somewhere} in the subgraph.
In fact, we can find most of these cycles even if we fix a vertex
$v$ and insist that a cycle of desired length passes through this
vertex.

Theorem \ref{pancyclic_thm4} says that for any fixed integer $l \geq
3$, if $p \gg n^{-1+1/(l-1)}$ then $G(n,p)$ almost surely has local
resilience $(1/2 + o(1))np$ with respect to containing cycles of
length $t$ for $l \leq t \leq n$. By carefully examining the proof,
one can realize that when finding middle length and long cycles, we
can insist on the cycle to pass through a fixed vertex. Thus we have
that for any fixed vertex $v$, there is a cycle of length $t$ for
$2l-1 \leq t \leq n$ which passes through $v$. Observe in addition
that for every odd integer $l \leq t < 2l-1$, we cannot guarantee a
cycle of length $t$ through a fixed vertex $v$ as can be seen using
Lemma \ref{pancyclic_random_lemma10}. Let $G=G(n,p)$. This lemma
implies that for any fixed vertex $v$, $|N^{(i)}_G(v)| = o(n)$ for
all $1 \leq i \leq l-2$. Therefore typically the degrees inside
$N^{(i)}_G(v)$ are all $o(np)$, and we can delete every edge inside
these sets without violating the maximum degree condition on $H$ to
get a graph $G-H$ which does not contain a cycle of length $t$
through $v$.

Similarly, Theorem \ref{pancyclic_pseudo_thm01} says that if $k$ is either 3 or an even integer satisfying $k \geq 4$ and $G=(V,E)$ is a $(n,d,\lambda)$-graph satisfying $d^{k-1} /n \gg \lambda^{k-2}$, then $G$ has local resilience $(1/2 + o(1))d$ with respect to containing cycles of length $t$ for $k \leq t \leq n$. Again even if we fix a vertex we can force all the middle length and long cycles to contain it. Namely, for any fixed vertex $v$, there exists a cycle of length $k+1$ up to $n$ passing through $v$.

\subsection{Paths through a given pair of vertices}

Another possible and rather straightforward extension of our results
is to show that random and pseudo-random graphs are locally
resilient with respect to the following property: for every given
pair of vertices $u,v$ and for every given length $l\ge t$ (where
$t$ is a constant depending on our choice of parameters, quite
similarly to the situation in Theorems \ref{pancyclic_thm4} and
\ref{pancyclic_pseudo_thm01}) there is a path of length $l$ between
$u$ and $v$. We do not provide much details here, but here is a very
short sketch of the argument. For medium length paths (between $t$
and $\delta n$, for some constant $\delta>0$) the proof is obtained
by a rather trivial modification of the corresponding proofs for
medium length cycles in Theorems \ref{pancyclic_thm4} and
\ref{pancyclic_pseudo_thm01}. For example, in the pseudo-random
case, instead of growing sets $X_i(v)$, $Y_i(v)$ as in the proof of
Theorem \ref{pancyclic_pseudo_thm01}, we grow disjoint sets $X_i(u)$
and $X_i(v)$ till they reach substantial size, and then find a
subgraph $G_1$ with large minimum degree on at least $(1/2+\delta)n$
vertices in the neighborhood of $X_i(u)$. Since $|V(G_1)|\ge
(1/2+\delta)n$, the set $X_i(v)$ has a neighbor $w$ in $G_1$. Due to
expansion properties of $G_1$, there is a path of linear length in
it starting from $w$. This path can be used to find paths of medium
length between $u$ and $v$. For paths of linear length, the key is
the ability to find a Hamilton path through a given pair of vertices
$u,v$ in an edge deleted random or pseudo-random graph. Here we can
argue as follows. First, if $e=(u,v)\not\in E(G)\setminus E(H)$
(where $G$ is the original (pseudo-)random graph, and $H$ is the
graph of deleted edges meeting the imposed condition on maximum
degree), add $e$ to the graph; clearly nothing really changes in its
edge distribution. Then, find a path $P$ of linear length with $e$
in somewhere in the middle (i.e., some sizable distance from both
ends), and then grow $P$ and close it to a Hamilton cycle through
rotations and extensions as usually, each time forbidding to touch
an interval of constant length surrounding $P$; our expansion
assumptions enable easily to meet this restriction. The so obtained
Hamilton cycle $C$ is guaranteed to contain $e$. Finally, omit $e$
from $C$, thus getting a Hamilton path between $u$ and $v$.

\subsection{Open Problems}

We believe that Theorem \ref{pancyclic_pseudo_thm1} can be extended
(with appropriate adjustments) to cycles of an arbitrary but fixed
odd length. More specifically, it is plausible that for an odd $k\ge
5$, if $G$ is an $(n,d,\lambda)$-graph and $d^{k-1}n/\gg
\lambda^{k-2}$, then the local resilience of $G$ with respect to
containing a cycle of length $k$ is $(1/2-o(1))d$. The validity of
this conjecture would allow to extend the assertion of Theorem
\ref{pancyclic_pseudo_thm01} to all $k\ge 3$.

A more natural generalization of Theorem \ref{pancyclic_pseudo_thm1}
(actually, of its original global resilience form as in
\cite{MR2145507}) is the following conjecture.
\begin{conj} Let $k\geq 5$ be an odd integer and $G$ be a $(n,d,\lambda)$-graph satisfying $d^{k-1} / n \gg \lambda^{k-2}$. Then $G$ has global resilience $(1/4 + o(1))nd$ with respect to being $C_k$-free.
\end{conj}


\appendix

\section{Proof of Theorem \ref{pancyclic_pseudo_thm2}}

In this appendix, we illustrate how the proof Theorem
\ref{pancyclic_pseudo_thm2} follows from the results of Sudakov and
Vu \cite{MR2462249}. First we repeat the statement of the theorem
here.

\begin{thm} Fix $\epsilon, \epsilon'$ such that $0 \leq 5\epsilon' < \epsilon, k \geq 3$ and let $G$ be an $(n,\epsilon', d,\lambda)$-graph satisfying $d^{k-1}/n \gg \lambda^{k-2}$. Then for large enough $n$, $G$ has local resilience at least $(1/2 - \epsilon)d$ with respect to being hamiltonian.
\end{thm}

As mentioned above, this theorem is not part of the original paper. In fact, they proved the following theorem.

\begin{thm}
 Let $\epsilon>0$ be fixed and $G$ be a $(n, d, \epsilon)$-graph such that $d/\lambda > \log^2 n$. Then for large enough $n$, $G$ has local resilience at least $(1/2 - \epsilon)d$ with respect to being hamiltonian.
\end{thm}

Unfortunately we cannot directly apply this result as our graph is not regular and we don't have a bound on $d / \lambda$. But Sudakov and Vu proved this result as a corollary of the following two results which can be modified to work in our situation,

\begin{thm} \label{thm_appendix1} 
For any fixed $\epsilon > 0$ and sufficiently large $n$, the following holds. Let $G$ be a connected graph of order $n$ 
such that every subset $U$ of $G$ of size at most $n/\log^4 n$ satisfies $|N_G(U)| \geq \frac{\log^4 n}{15}\cdot |U|$ 
and every 
subset $W$ of size at least $n / \log^3 n$ has $|N_G(W)| \geq \frac{1+\epsilon}{2} n$. Then $G$ contains a Hamilton 
cycle.
\end{thm}
\begin{lemma} \label{lemma_appendix1} For any fixed $\epsilon > 0$ and sufficiently large $n$ the following holds. Let $G$ be a $(n,d,\lambda)$-graph with $d / \lambda > \log^2 n$ and let $H$ be a subgraph of $G$ with maximum degree at most $(1/2 - \epsilon)d$. Then the graph $G' = G-H$ is
\begin{itemize}
\item connected;
\item every subset $U$ of $G'$ of size at most $n/ \log^4 n$ satisfies $|N_{G'}(U)| \geq 
\frac{\log^4 n}{15}\cdot |U|$.
\item every subset $W$ of size at least $n/ \log^3 n$ has $|N_{G'}(W)| \geq \frac{1+\epsilon}{2}n$.
\end{itemize}
\end{lemma}

By applying Lemma \ref{pancyclic_pseudo_lemma10} and Lemma \ref{pancyclic_pseudo_lemma1} we have the following result similar to Lemma \ref{lemma_appendix1}.
\begin{lemma}\label{lemma_appendix2} Fix $\epsilon, \epsilon'$ such that $0 \leq 5\epsilon' < \epsilon, k \geq 3$ and let $G$ be an $(n,\epsilon', d,\lambda)$-graph satisfying $d^{k-1} / \lambda^{k-2} \gg n$ and let $H$ be a subgraph of $G$ with maximum degree at most $(1/2 - \epsilon)d$. Then the graph $G' = G-H$ is
\begin{itemize}
\item connected;
\item every subset $U$ of $G'$ of size at most $(\lambda^2 / d^2) n$ satisfies $|N_{G'}(U)| \geq \frac{\epsilon d^2}{4\lambda^2}\cdot |U|$.
\item for an arbitrary function $\delta(n)$ growing to infinity, every subset $W$ of size at least $\delta(n) (\lambda^2 / d^2) n$ has $|N_{G'}(W)| \geq \frac{1+\epsilon}{2}n$.
\end{itemize}
\end{lemma}

Therefore we only need to prove the following theorem which is a variant of Theorem \ref{thm_appendix1}.
\begin{thm} Let $\epsilon > 0, k \geq 3$ be fixed $d^{k-1}/\lambda^{k-2} \gg n$. Then for sufficiently large $n$ the following holds. Let $G$ be a connected graph of order $n$
such that every subset $U$ of $G$ of size at most $(\lambda^2/d^2)n$
satisfies $|N_G(U)| \geq \frac{\epsilon d^2}{4\lambda^2}\cdot |U|$
and every subset $W$ of size at least $\delta(n)(\lambda^2/d^2) n$
has $|N_G(W)| \geq \frac{1+\epsilon}{2} n$. Then $G$ contains a
Hamilton cycle.
\end{thm}
We omit the proof which is a word by word translation of the proof of Theorem \ref{thm_appendix1}. (Actually we are in a more simple situation since we only need $k/2$ rotations compared to $\log n / \log \log n$ as in the original proof.)

\bibliographystyle{plain}
\bibliographystyle{plain}

\begin{thebibliography}{99}

\bibitem{MR1302331}
N. Alon,
\newblock Explicit {R}amsey graphs and orthonormal labelings,
\newblock {\em Electron. J. Combin.} 1 (1994), Research Paper 12, 8 pp.

\bibitem{MR000000}
S. Ben-Shimon, M. Krivelevich and B. Sudakov,
\newblock {\em in preparation}.

\bibitem{MR0285424}
J.A. Bondy,
\newblock Pancyclic graphs {I},
\newblock {\em J. Combinatorial Theory Ser. B} 11 (1971), 80--84.

\bibitem{MR0340095}
J.A. Bondy and M. Simonovits,
\newblock Cycles of even length in graphs,
\newblock {\em J. Combinatorial Theory Ser. B} 16 (1974), 97--105.

\bibitem{MR1611825}
S. Brandt, R. Faudree and W. Goddard,
\newblock Weakly pancyclic graphs,
\newblock {\em J. Graph Theory} 27 (1998), 141--176.

\bibitem{MR2383452}
D. Dellamonica, Y. Kohayakawa, M. Marciniszyn and A. Steger,
\newblock On the resilience of long cycles in random graphs,
\newblock {\em Electron. J. Combin.} 15 (2008), Research Paper 32.

\bibitem{MR2159259}
R. Diestel,
\newblock {\bf Graph theory}, Volume 173 of {\em Graduate Texts in
  Mathematics},
\newblock Springer-Verlag, Berlin, 3rd edition, 2005.

\bibitem{MR2430433}
A. Frieze and M. Krivelevich,
\newblock On two {H}amilton cycle problems in random graphs,
\newblock {\em Israel J. Math.} 166 (2008), 221--234.

\bibitem{MR1339852}
P.E. Haxell, Y. Kohayakawa and T.{\L}uczak,
\newblock Tur\'an's extremal problem in random graphs: forbidding even
cycles,
\newblock {\em J. Combin. Theory Ser. B} 64 ( 1995), 273--287.

\bibitem{MR1394514}
P.E. Haxell, Y. Kohayakawa and T. {\L}uczak.
\newblock Tur\'an's extremal problem in random graphs: forbidding odd
cycles,
\newblock {\em Combinatorica} 16 (1996), 107--122.

\bibitem{MR1782847}
S. Janson, T. {\L}uczak and A. Rucinski,
\newblock {\bf Random graphs},
\newblock Wiley-Interscience Series in Discrete Mathematics and Optimization.
  Wiley-Interscience, New York, 2000.

\bibitem{MR1945368}
J. Kim, B. Sudakov and V. Vu,
\newblock On the asymmetry of random regular graphs and random
graphs,
\newblock {\em Random Structures and Algorithms} 21 (2002), 216--224.

\bibitem{MR2223394}
M. Krivelevich and B. Sudakov,
\newblock Pseudo-random graphs,
\newblock In: {\em More sets, graphs and numbers}, Volume~15 of {\em Bolyai Soc.
  Math. Stud.}, pages 199--262. Springer, Berlin, 2006.

\bibitem{MR2321240}
L. Lov{\'a}sz,
\newblock {\bf Combinatorial problems and exercises},
\newblock AMS Chelsea Publishing, Providence, RI, 2nd edition, 2007.

\bibitem{MR2145507}
B. Sudakov, T. Szab{\'o} and V. Vu,
\newblock A generalization of {T}ur\'an's theorem,
\newblock {\em J. Graph Theory} 49 (2005), 187--195.

\bibitem{MR2462249}
B. Sudakov and V. Vu,
\newblock Local resilience of graphs,
\newblock {\em Random Structures and Algorithms} 33 (2008), 409--433.

\bibitem{MR0018405}
P. Tur{\'a}n,
\newblock Eine {E}xtremalaufgabe aus der {G}raphentheorie,
\newblock {\em Mat. Fiz. Lapok} 48 (1941), 436--452.

\bibitem{MR0318000}
D. Woodall,
\newblock Sufficient conditions for circuits in graphs,
\newblock {\em Proc. London Math. Soc. } 24 (1972), 739--755.

\end{thebibliography}

\end{document}